\documentclass[a4paper,draft]{article}
\usepackage{amsmath}
\usepackage{amsthm}
\usepackage{amssymb}
\usepackage{color}

\newtheorem{theorem}{Theorem}[]
\newtheorem{lemma}{Lemma}[]
\newtheorem{problem}{Problem}[]
\newtheorem*{kn_problem}{Known problem}
\newtheorem*{grand}{}

\mathversion{bold}
\begin{document}
\title{The concrete theory of numbers: initial numbers and wonderful properties of numbers repunit}
\author{Tarasov, B.\,V.\thanks{Tarasov, B.\,V.
The concrete theory of numbers: initial numbers and wonderful properties of numbers repunit.
MSC 11A67+11B99.
\textcopyright 2007 Tarasov, B.\,V.,
independent researcher.
}\\
}
\maketitle
\begin{abstract}
In this work initial numbers and repunit numbers have been studied. All numbers have been
considered in a decimal notation. The problem of simplicity of initial numbers has been studied.
Interesting properties of numbers repunit are proved: $gcd(R_a, R_b) = R_{gcd(a,b)}$;
 $R_{ab}/(R_aR_b)$ is an integer only if $gcd(a,b) = 1$, where
 $a\geq1$, $b\geq1$~are integers. Dividers of numbers repunit, are researched by a degree of prime number.
\end{abstract}
\begin{flushright}
\textbf{\textcolor[rgb]{1.00,0.00,0.50}{Devoted to the tercentenary from the date of birth (4/15/1707)
\\of Leonhard Euler \\}}
\end{flushright}

{\centering\section[]{Introduction}\par}

Let $x\geq0$, $n\geq0$~be integers.  An integer $N$, which record consists from $n$ \emph{records}
of number $x$, we shall designate by
\begin{equation}N=\{x\}_n=x\ldots x,~n>0.\end{equation}
For $n=0$ it is received $\{x\}_0=\varnothing$ an empty record. For example, $\{10\}_3 1=1010101$, $\{10\}_0 1=1$, etc.\par
Palindromic numbers of a kind
\begin{equation}E_{n,k}=\{1\{0\}_k\}_n1,\end{equation}
where $n\geq0$, $k\geq0$ we will name initial numbers. We will notice that $E_{0,k}=1$
at any $k\geq0$.\par
Numbers repunit(see\cite{Graham,Weisstein,The Prime Clossary repunit}) are natural numbers, which records consist
of units only, i.e. by definition
\begin{equation}R_n=E_{n-1,0},\end{equation}
where $n\geq1$.
\par In decimal notation the general formula for numbers repunit is
\begin{equation}R_n=(10^n-1)/9,\end{equation}
where $n=1,2,3,\ldots$~.
\par There are known only five prime repunit for $n=$2,19, 23, 317, 1031.
\par \begin{kn_problem}[(Prime repunit numbers\cite{Weisstein})]~Whether exists infinite number
of prime numbers repunit ?\end{kn_problem}
\par Will we use designations further :
\par $(a,b)=gcd(a,b)$~the greatest common divider of integers $a>0$, $b>0$.
\par $p,q$ odd prime numbers.
\par If it is not stipulated specially, the integer positive numbers are considered.\par

{\centering\section{Initial numbers}}\par

Let's consider the trivial properties of initial numbers.\par
\begin{theorem}[]Following trivial statements are fair :\par
(1) General formula of initial numbers is
\begin{equation}E_{n,k}=\frac{R_{(k+1)(n+1)}}{R_{k+1}}=\frac{10^{(k+1)(n+1)}-1}{10^{k+1}-1}.\end{equation}\par
(2) For $k\geq0$, $n\geq m\geq1$ if $n+1\equiv0(mod\,(m+1))$,\\
 then $(E_{n,k},E_{m,k})=E_{m,k}$.\par
(3) For $k\geq0$, $n>m\geq1$ if integer $s\geq1$, exists such\\ that $n+1\equiv0(mod\,(s+1))$,
 $m+1\equiv0(mod\,(s+1))$, then\\ $(E_{n,k},E_{m,k})\geq E_{s,k}>1$.\par
(4) For $k\geq0$, $n>m\geq1$ $(E_{n,k},E_{m,k})=1$  when and only then,\\ $(n+1,m+1)=1$.
\end{theorem}
\begin{proof}[\textbf{Proof}]
1) Properties (1)---(3) are obvious.\par
2) The Proof of property (4). \textbf{Necessity}. Let \\$(E_{n,k},E_{m,k})=~1$ and $(n+1,m+1)=s>1$, $s-1\geq1$.
From property (3) of the theorem follows that $(E_{n,k},E_{m,k})\geq E_{s-1,k}=\{1\{0\}_k\}_{s-1}1>1$. Appears the contradiction
.\par
\textbf{Sufficiency} of property (4). Let $(n+1,m+1)=1$, then will be integers $a>0$, $b>0$, such that either
$a(n+1)=b(m+1)+1$ \\or $b(m+1)=a(n+1)+1$. Let's assume, that $(E_{n,k},E_{m,k})=d>1$.\par
a) Let $a(n+1)=b(m+1)+1$, then $E_{b(m+1),k}=E_{a(n+1)-1,k}=(10^{a(n+1)(k+1)}-1)/(10^{k+1}-1)\equiv0(mod\, E_{n,k})
\equiv0(mod\, d)$.\par
On the other hand $E_{b(m+1),k}=(10^{(k+1)\{b(m+1)+1\}}-1)/(10^{k+1}-1)=((10^{b(m+1)(k+1)}-1)/(10^{k+1}-1))
\cdot10^{k+1} + 1\equiv$ \\ $\equiv1(mod\, E_{m,k})\equiv1(mod\, d)$. Appears the contradiction.\par
b) Let $b(m+1)=a(n+1)+1$, then $E_{a(n+1),k}=E_{b(m+1)-1,k}=(10^{b(m+1)(k+1)}-1)/(10^{k+1}-1)\equiv0(mod\, E_{m,k})
\equiv0(mod\, d)$.\par
On the other hand $E_{a(n+1),k}=(10^{(k+1)\{a(n+1)+1\}}-1)/(10^{k+1}-1)=((10^{a(n+1)(k+1)}-1)/(10^{k+1}-1))
\cdot10^{k+1} + 1\equiv$ \\ $\equiv1(mod\, E_{n,k})\equiv1(mod\, d)$. Have received the contradiction. \end{proof}

{\centering\section{Numbers repunit}}\par

Let's consider trivial properties of numbers repunit.\par
\begin{theorem}[]Following trivial statements are fair :\par
(1) The number $R_n$ is prime only if $n$ number is prime.\par
(2) If $p>3$ all prime dividers of number $R_p$ look like $1+2px$ where $x\geq1$ is integer.\par
(3) $(R_a,R_b)=1$ if and only if $(a,b)=1$.
\end{theorem}
\begin{proof}[\textbf{Proof}]
Property (1) of theorem is proved in (\cite{Graham,Weisstein}), property (2) is proved in (\cite{Vinogradov}), as exercise.
Property (3) is the corollary of the theorem 1.
\end{proof}

\begin{theorem}[] $(R_a,R_b)=R_{(a,b)}$, where $a\geq1$, $b\geq1$ are integers.\par
\end{theorem}
\begin{proof}[\textbf{Proof}]
Validity of the theorem for $(a,b)=1$ follows from property (3) of theorem\,2. Let $(a,b)=d>1$, where $a=a_1d$, $b=b_1d$,
$(a_1,b_1)=1$. Let's consider equations
\begin{displaymath}R_a=R_d\cdot\{10^{d(a_1-1)}+\,\ldots\,+10^d+1\},\end{displaymath}
\begin{displaymath}R_b=R_d\cdot\{10^{d(b_1-1)}+\,\ldots\,+10^d+1\}.\end{displaymath}
\par Let
\begin{displaymath}A=10^{d(a_1-1)}+\,\ldots\,+10^d+1,\end{displaymath}
\begin{displaymath}B=10^{d(b_1-1)}+\,\ldots\,+10^d+1.\end{displaymath}
\par Let's assume, that $(A,B)>1$, and $q$ is a prime odd number such that
\begin{equation}A\equiv0(mod\, q), ~B\equiv0(mod\, q).\end{equation}
\par If $q=3$, then $10^t\equiv1(mod\, q)$ for any integer $t\geq1$. Then from (6) it follows that
$a_1\equiv b_1\equiv0(mod\, q)$. Have received the contradiction.\par
Thus, $q>3$. Then there exists an index $d_{min}$, to which the number $10^d$ belongs on the module $q$.
\begin{displaymath}(10^d)^{d_{min}}\equiv1(mod\, q),\end{displaymath}
where $d_{min}\geq1$.
\par If $d_{min}=1$, then it follows from (6) that $a_1\equiv b_1\equiv0(mod\, q)$. Have received the contradiction. Hence
$d_{min}>1$. As $R_a\equiv R_b\equiv~0(mod\, q)$, then $(10^d)^{a_1}\equiv1(mod\, q)$ and
$(10^d)^{b_1}\equiv1(mod\, q)$.\\ Then $a_1\equiv b_1\equiv0(mod\, d_{min})$. Have received the contradiction.
\end{proof}

\begin{theorem}[] Let $p>3$ be a prime number, $k\geq t\geq1$, $t\geq s\geq1$ integer numbers. Then
\begin{equation}\label{e:T4}gcd(R_{p^k}/R_{p^t}, R_{p^s})=1.\end{equation}
\end{theorem}
\begin{proof}[\textbf{Proof}]
Let's consider expression
\begin{displaymath}A=R_{p^k}/R_{p^t}=(10^{p^t})^{p^{k-t}-1} + (10^{p^t})^{p^{k-t}-2} +\,\ldots\,+ 10^{p^t} + 1.
\end{displaymath}
\par If $(A,R_{p^s})>1$, then the prime number $q$ exists such that\par
 $A\equiv~0(mod\, q)$ è $R_{p^s}\equiv~0(mod\, q)$.
 Hence $10^{p^t}\equiv~1(mod\, q)$, then $A\equiv~p^{k-t}\equiv~0(mod\, q)$, $p=q=3$.
Have received the contradiction, because $p>3$.
\end{proof}
\begin{theorem}[] Let $a\geq1$, $b\geq1$ are integers, then the following statements are true :\par
(1) If $(a,b)=1$, then
\begin{equation}gcd(R_{ab},R_aR_b)=R_aR_b.\end{equation}
\par (2) If $(a,b)>1$, then
\begin{equation}R_aR_b /R_{(a,b)} \leq gcd(R_{ab},R_aR_b)<R_aR_b.\end{equation}

\end{theorem}
\begin{proof}[\textbf{Proof}]
1) Let $(a,b)=1$, then $(R_a,R_b)=R_{(a,b)}=1$,\\ $R_{ab}=R_aX=R_bY$, $X=cR_b$, where $c\geq1$ is integer.
$R_{ab}=cR_aR_b$.\par
2) Let $(a,b)=d>1$, $a=a_1d$, $b=b_1d$, $(a_1,b_1)=1$, $a_1\geq1$, $b_1\geq1$.
As $gcd(R_a,R_b)=R_{(a,b)}$, we receive equality
\begin{equation}R_a=R_{(a,b)}X, R_b=R_{(a,b)}Y,\end{equation}
where $(X,Y)=1$.\par
Further, $R_{ab}=R_aA=R_bB=XAR_{(a,b)}=YBR_{(a,b)}$, $XA=YB$,\\
$A=Yz$, $B=Xz$, $z\geq1$ is integer. Then $R_{ab}=XYR_{(a,b)}z$,\\
$R_{ab}=zR_aR_b /R_{(a,b)}$. We have proved, that\\ $R_aR_b /R_{(a,b)} \leq gcd(R_{ab},R_aR_b)$.
\par Let's assume, that $gcd(R_{ab},R_aR_b)=R_aR_b$, then $R_{ab}=zR_aR_b$, where $z\geq1$  is integer.
Let's consider equalities
\begin{displaymath}R_{ab}=R_aA=R_bB,\end{displaymath}
\par where
\begin{displaymath}A=10^{a(b-1)}+10^{a(b-2)}+\,\ldots\,+10^a+1,\end{displaymath}
\begin{displaymath}B=10^{b(a-1)}+10^{b(a-1)}+\,\ldots\,+10^b+1.\end{displaymath}
Since $A=R_bz$, $B=R_az$, $10^a\equiv1(mod\,R_{(a,b)})$,\\ $10^b\equiv1(mod\,R_{(a,b)})$,
then $A\equiv B\equiv0(mod\,R_{(a,b)})$, hence \\ $a\equiv b\equiv0(mod\,R_{(a,b)})$.\par
Thus, comparison $(a,b)\equiv0(mod\,R_{(a,b)})$ or\\ $d\equiv0(mod\,R_{d})$ is fair,
that contradicts an obvious inequality
\begin{equation}(10^x-1)/9 >x,\end{equation}
where $x>1$ is real. \end{proof}
\begin{grand}[\{$\bigstar$\} \textbf{The Important corollary of the theorem 5}] \ \\ Number $R_{ab}/(R_aR_b)$ is integer when and only when
$(a,b)=1$, where $a\geq1$, $b\geq1$ are integers.\end{grand}
Let's quote some trivial statements for numbers repunit.
\begin{lemma}[]If $a=3^nb$, $(b,3)=1$, then
\begin{equation}\label{e:L1}R_a\equiv0(mod\,3^n),\,but\,R_a\not\equiv0(mod\,3^{(n+1)}).\end{equation}
\end{lemma}\par
\begin{proof}[\textbf{Proof}]\par
If $n=1$, then $R_a=R_3B$, where $B=10^{3(b-1)}+\,\ldots\,+10^3+1$, $R_3\equiv0(mod\,3)$, $B\equiv b\not\equiv0(mod\,3)$.
Thus,\\ $R_a\equiv0(mod\,3)$, but $R_a\not\equiv0(mod\,3^2)$.
\par Let comparisons \eqref{e:L1} be proved for $n\leq k-1$. We shall consider $a=3^kb$, $(b,3)=1$. Then $R_a=R_{3^{k-1}b}A$, where
$A=10^{3^{k-1}b2}+10^{3^{k-1}b}+1$. \par$R_{3^{k-1}b}\equiv0(mod\,3^{k-1})$, but $R_{3^{k-1}b}\not\equiv0(mod\,3^{k})$,
 $A\equiv0(mod\,3)$, but $A\not\equiv0(mod\,3^2)$.
\end{proof}
\begin{lemma}[] If $n\geq0$ is integer, then
\begin{equation}\label{e:L2}r_n=10^{11^n}+1\equiv0(mod\,11^{n+1}),\,but\,r_n\not\equiv0(mod\,11^{n+2}).\end{equation}
\end{lemma}\par
\begin{proof}[\textbf{Proof}]
$r_0=11\equiv0(mod\,11)$, but $r_0=11\not\equiv0(mod\,11^2)$.\\
$r_1=10^{11}+1\equiv0(mod\,11^2)$, but $r_1\not\equiv0(mod\,11^3)$.
\par Let's make the inductive assumption, that formulas \eqref{e:L2} are proved for\\ $n\leq k-1$, where $k-1\geq1$, $k\geq2$.
Let $n=k$, then\\ $r_k=10^{11^k}+1=(10^{11^{k-1}})^{11}+1=r_{k-1}A$, where
\begin{multline}\label{e:L21}A=10^{11^{k-1}10}-10^{11^{k-1}9}+10^{11^{k-1}8}-10^{11^{k-1}7}
+10^{11^{k-1}6}-\\
-10^{11^{k-1}5}+10^{11^{k-1}4}-10^{11^{k-1}3}+10^{11^{k-1}2}-10^{11^{k-1}}+1.\end{multline}
Since, due to the inductive assumption $10^{11^{k-1}}\equiv-1(mod\,11^k)$, where\\ $k\geq2$, then $A\equiv11(mod\,11^k)$.
Then $A\equiv0(mod\,11)$, but\\ $A\not\equiv0(mod\,11^2)$.
Thus, we receive, that $r_k\equiv0(mod\,11^{k+1})$, but\\ $r_k\not\equiv0(mod\,11^{k+2})$.
\end{proof}
\begin{lemma}[] For an integer $a\geq1$, the following statements are true :\par
(1) If $a$ is odd, then $R_a\not\equiv0(mod\,11)$.\par
(2) If $a=2(11^n)b$, $(b,11)=1$, then
\begin{equation}\label{e:L3}R_a\equiv0(mod\,11^{n+1}),\,but\,R_a\not\equiv0(mod\,11^{n+2}).\end{equation}
\end{lemma}\par
\begin{proof}[\textbf{Proof}]
If $a$ is odd, then $R_a\equiv1(mod\,11)$.
If $a=2(11^n)b$, $(b,11)=1$, then
$R_a=((10^{2(11)^n})^b-1)/9=R_{11^n}\cdot r_n\cdot A$,
where $r_n=10^{11^n}+1$, \\$A=10^{2(11^n)(b-1)}+\,\ldots\,+10^{2(11^n)}+1$. $R_{11^n}\not\equiv0(mod\,11)$,\\
$A\equiv ~b\not\equiv~0(mod\,11)$. Then validity of the statement (2) of lemma 3 follows from lemma 2.
\end{proof}
\begin{grand}[\{$\bigstar$\} \textbf{The assumption: }\textbf{the general formula for} $gcd(R_{ab},R_aR_b)$] \ \\
If $a\geq1$, $b\geq1$ are integers, $d=(a,b)$, where $d=3^L\cdot11^S\cdot c$, $(c,3)=1$, $(c,11)=1$, $L\geq0$,
$S\geq0$,
then equalities are true :\par
--- if $c$ is an odd number, then
\begin{equation}\label{e:guess1}gcd(R_{ab},R_aR_b)=((R_aR_b)/R_{(a,b)})\cdot3^L,\end{equation}
\par--- if $c$ is an even number, then
\begin{equation}\label{e:guess2}gcd(R_{ab},R_aR_b)=((R_aR_b)/R_{(a,b)})\cdot3^L\cdot11^S.\end{equation}
\end{grand}
Let's give another two obvious statements in which divisors of numbers\\ repunit are studied, as degrees of prime number.
\begin{lemma}[] If $p,q$ are prime numbers and $R_p\equiv0(mod\,q)$, but \\$R_p\not\equiv0(mod\,q^2)$, then statements are true :\par
(1) For any integer $r$, $0<r<q$, $R_{pr}\not\equiv0(mod\,q^2)$.\par
(2) For any integer $n$, $n\geq1$, $R_{p^n}\not\equiv0(mod\,q^2)$.
\end{lemma}\par
\begin{proof}[\textbf{Proof}]
1) $R_{pr}=R_p\cdot\widehat{R}_{pr}$, where $\widehat{R}_{pr}=10^{p(r-1)}+10^{p(r-2)}+\,\\+\,\ldots\,+10^{p}+1$.
If $R_{pr}\equiv0(mod\,q^2)$, then $\widehat{R}_{pr}\equiv0(mod\,q)$,\\ $r\equiv~0(mod\,q)$. Have received the contradiction.\par
2) If $n>1$ found such that $R_{p^n}\equiv0(mod\,q^2)$, then from \eqref{e:T4} follows $(R_{p^n}/R_p, R_p)=1$. Have received the contradiction.
\end{proof}
\begin{lemma}[] If $p,q$ are prime numbers and $R_p\equiv0(mod\,q)$, then\par $R_{pq^n}\equiv0(mod\,q^{n+1})$.
\end{lemma}\par
\begin{proof}[\textbf{Proof}]
Since $R_{pq}=R_p\cdot\widehat{R}_{pq}$, where $\widehat{R}_{pq}=10^{p(q-1)}+\\+\,10^{p(q-2)}+\,\ldots\,+10^{p}+1$, then
$\widehat{R}_{pq}\equiv0(mod\,q)$, $R_{pq}\equiv0(mod\,q^2)$.\\
Let's assume that $R_{pq^{n-1}}\equiv0(mod\,q^{n})$. Then\\
$R_{pq^n}=R_{pq^{n-1}\cdot q}=R_{pq^{n-1}}\cdot\widehat{R}_{pq^{n-1}\cdot q}$,
where\\ $\widehat{R}_{pq^{n-1}\cdot q}=10^{pq^{n-1}\cdot(q-1)}+10^{pq^{n-1}\cdot(q-2)}+\,\ldots\,+10^{pq^{n-1}}+1\equiv~0(mod\,q)$,\\
$R_{pq^{n}}\equiv0(mod\,q^{n+1})$.
\end{proof}

{\centering\section{Problem of simplicity of initial numbers}}\par
Let's consider the problem of simplicity of initial numbers $E_{n,k}$, where\\ $k\geq0$, $n\geq0$.\par
If $k=0$, then $E_{n,0}=R_{n+1}$. Thus, simplicity of numbers $E_{n,0}$ -- is known problem of prime numbers repunit $R_p$,
where $p$ is prime number.\par
If $n=1$, then $E_{1,k}=1\{0\}_k1=10^{k+1}+1$. As number $E_{1,k}$ can be prime only when $k+1=2^m$, $m\geq0$ is integer,
 then we come to the known problem of simplicity of the generalized Fermat numbers $f_m(a)=a^{2^m}+1$ for $a=10$.
 Generalized Fermat numbers nave been define by Ribenboim
 \cite{Ribenboim} in 1996, as numbers of the form $f_n(a)=a^{2^n}+1$, where $a>2$ is even.\\
The generalized Fermat numbers $f_n(10)=10^{2^n}+1$ for $n\leq14$ are prime only if $n=0,1$. $f_0(10)=11$, $f_1(10)=101$.\par
\begin{theorem}[]Let $n>1$, $k>0$. If any of conditions \par
(1) $n$ number is odd, \par
(2) $k$ number is odd, \par
(3) $n+1\equiv0(mod\,3)$, \par
(4) $(n+1,k+1)=1$, \\
is true, then number $E_{n,k}$ is compound.
\end{theorem}
\begin{proof}[\textbf{Proof}]
1) $n+1=2t$, $t>1$. Then $E_{n,k}=E_{t-1,k}\cdot(10^{t(k+1)}+1)$, where $t>1$, $t-1\geq1$. As $E_{t-1,k}>1$, then $E_{n,k}$ is compound number.\par
2) Let $k$ be an odd number. Due to the proved condition (1) we count that number $(n+1)$ is odd. $k+1=2t\geq2$, $t\geq1$. Further,
\begin{displaymath}E_{n,k}=E_{n,t-1}\cdot((10^{(n+1)t}+1)/(10^t+1)),\end{displaymath}
where $n>1$, $t-1\geq0$, $E_{n,t-1}>1$, number $(10^{(n+1)t}+1)/(10^t+1)>1$ is integer.\par
3) If $n+1\equiv0(mod\,3)$, then $E_{n,k}\equiv0(mod\,3)$, $E_{n,k}>11$.\par
4) Let $n>1$, $k\geq1$, $(n+1,k+1)=1$, then
\begin{displaymath}E_{n,k}=R_{(n+1)(k+1)}/R_{(k+1)}=R_{(n+1)}\cdot(R_{(n+1)(k+1)}/(R_{k+1}\cdot R_{n+1})).\end{displaymath}
Due to the theorem 5 number $z=R_{(n+1)(k+1)}/(R_{k+1}\cdot R_{n+1})$ is integer. Further, \\
$z>(10^{(n+1)(k+1)}-1)/(10^{n+k+2})=10^{nk-1} - 1/(10^{n+k+2})$, $nk-1\geq1$, thus, $z>1$.
\end{proof}
Question of simplicity of initial numbers under conditions, when\\
$(n+1,k+1)>1$, $(n+1)$ number is odd, $(k+1)$ number is odd,\\ $n+1\not\equiv0(mod\,3)$, remains open.\par
In particular, it is interesting to considerate numbers $E_{p-1,p-1}=R_{p^2}/R_p$, where $p$ is prime number. For $p<100$ numbers $E_{p-1,p-1}$ are compound.

{\centering\section{The open problems of numbers repunit}}\par

The known problem of numbers repunit remains open.
\begin{problem}[(Prime repunit numbers\cite{Weisstein})]~Whether there exists infinite number of prime numbers
$R_p$, $p$--prime number ?\end{problem}
\begin{problem}[]~Whether all numbers $R_p$, $p$--prime number, are numbers free from squares ?\end{problem}
The author has checked up for $p<97$, that numbers $R_p$ are free from squares.
Another following open questions are interesting :
\begin{problem}[]~If number $R_p$ is free from squares, where $p>3$ is prime number, whether will number $n$, be found such what
number $R_{p^n}$ contains a square~?\end{problem}
\begin{problem}[]~$p$ is prime number, whether there are simple numbers of a kind\\
$E_{p-1,p-1}=R_{p^2}/R_p$ ?\end{problem}
The author has checked up to $p\leq127$, that numbers $E_{p-1,p-1}$ is compound.
It is known, that $R_p$ divide by  number $(2p+1)$ for prime numbers $p=41,53$, $R_p$ divide by number $(4p+1)$ for prime numbers $p=13,43,79$.\\
There appears a question\,:
\begin{problem}[]~Whether there is infinite number of prime numbers $p$, such that $R_p$ divide by number $(2p+1)$ or is number
$(4p+1)$ ?\end{problem}
\begin{grand}[\textbf{The remark}]~If the number $p>5$ Sophie Germain prime (i.e. number $2p+1$ is prime too), then either $R_p$
or $R_{p}^{+}=(10^p+1)/11$ divide by number $(2p+1)$.\end{grand}

{\centering\section{The conclusion}}\par

\textbf{\textcolor[rgb]{1.00,0.00,0.00}{Leonhard Euler}}, professor of the Russian Academy of sciences since 1731,
\textbf{\textcolor[rgb]{1.00,0.00,0.00}{has paid mathematics forever !}}
Euler's invisible hand directs the development of concrete mathematics for more than 200 years. \par
Euler's titanic work which has opened a way to freedom to mathematical community, admires.
The pleasure caused by Euler's works warms hearts.
\pagebreak

\par

\textcolor[rgb]{0.00,0.00,1.00}{---------------------------------------------------------------------}\par
\textcolor[rgb]{0.00,0.00,1.00}{Institute of Thermophysics, Siberian Branch of RAS }\par
\textcolor[rgb]{0.00,0.00,1.00}{Lavrentyev Ave., 1, Novosibirsk, 630090, Russia }\par
\textcolor[rgb]{0.00,0.00,1.00}{E-mail: tarasov@itp.nsc.ru }\par
\textcolor[rgb]{0.00,0.00,1.00}{---------------------------------------------------------------------}\par
\textcolor[rgb]{0.00,0.00,1.00}{Independent researcher, }\par
\textcolor[rgb]{0.00,0.00,1.00}{E-mail: tarasov-b@mail.ru }\par
\textcolor[rgb]{0.00,0.00,1.00}{---------------------------------------------------------------------}\par


\begin{thebibliography}{9}
\par
\bibitem{Vinogradov} Vinogradov I.\,M.
\textit{Osnovy teorii chisel.\,-\,M.\,:\,Nauka,\,1981.}
\bibitem{Graham} Ronald L.\,Graham,\,Donald E.\,Knuth,\,Oren Patashnik,\\
\textit{Concrete Mathematics :\,A Foundation for Computer Science,\,2nd edition (Reading,\,Massachusetts:\,Addison-Wesley), 1994.}
\bibitem{Weisstein} Weisstein,\,Eric W.
\textit{"Repunit." From MathWorld--A Wolfram Web Resource.
---http://mathworld.wolfram.com/Repunit.html/.\\
\textcopyright 1999---2007 Wolfram Research,\,Inc.}
\bibitem{The Prime Clossary repunit}
\textit{ The Prime Clossary repunit.\\
---http://primes.utm.edu/glossary/page.php?sort=Repunit/.}
\bibitem{Ribenboim} Ribenboim,\,P.
\textit{"Fermat Numbers"\,and\,"Numbers $k\times2^n\pm1$." §2.6 and 5.7 in The New Book of Prime Number Records.
New York: Springer-Verlag, pp. 83-90 and 355-360, 1996.}
\end{thebibliography}
\end{document}